\documentclass[12pt]{article}
\usepackage{amsfonts}
\usepackage{amsmath}
\usepackage{amssymb}
\usepackage{amsthm}
\usepackage[english,russian]{babel}
\usepackage[T2A]{fontenc}

\newtheorem {thm}{Теорема}
\newtheorem {lem}{Лемма}
\newtheorem {sta}{Утверждение}

\theoremstyle{definition}
\newtheorem {rem}{Замечание}
\newtheorem {defin}{Определение}

\DeclareMathOperator{\aff}{\mathit {aff}}

\title{Критерий приводимости параллелоэдра}
\author{А.~Ордин, А.~Магазинов\thanks{Второй автор поддержан Правительством РФ, проект 11.G34.31.0053 и грантом РФФИ 11-01-00633.}}

\begin{document}

\maketitle

\subsection*{Аннотация}

Параллелоэдр называется {\it приводимым}, если его можно представить в виде прямого произведения двух параллелоэдров меньшей размерности. Критерий приводимости параллелоэдра в терминах графа Венкова был доказан первым автором в 2005 году в его диссертации. В несколько переработанном виде оригинальное доказательство было
изложено вторым автором на семинаре ``Дискретная геометрия и геометрия чисел'' (МГУ) в декабре 2011 года. Настоящая работа является развернутым конспектом
упомянутого доклада.

\subsection*{Abstract}

\subsubsection*{A criterion of reducibility for a parallelohedron\newline {\rm A.~Ordine, A.~Magazinov}}
A parallelohedron is called {\it reducible}, if it can be represented as a direct product of two parallelohedra of lower dimension. In his Ph.D. thesis (2005) 
the first author proved a criterion of reducibility of a parallelohedron in terms of the Venkov graph. In December 2011 the second author presented a
slightly revised version of the original proof at the seminar ``Discrete Geometry and Geometry of Numbers'' (Moscow State University). The present
paper follows that talk.

\section{Введение}

Основными объектами, изучаемыми в настоящей статье, являются параллелоэдры.

\begin{defin}
{\it Параллелоэдром} (см.~\cite{Fed1885}) размерности $d$ называется такой многогранник $P$, что существует разбиение грань-в-грань пространства 
$\mathbb R^d$ транслятами $P$.
\end{defin}

Сформулируем условия для выпуклого многогранника $P$, необходимые и достаточные для того, чтобы $P$ являлся параллелоэдром. Для этого потребуется
Определение \ref{def:2}.

\begin{defin}\label{def:2}
Пусть $Q$ --- многогранник в $\mathbb R^d$, все гиперграни которого имеют центр симметрии, $F$ --- любая $(d-2)$-мерная грань $Q$. Тогда множество
всех гиперграней $Q$, параллельных $F$, называется {\it пояском} $Q$.
\end{defin}

\begin{rem}
В силу центральной симметрии, любая гипергрань пояска содержит ровно две $(d-2)$-мерные грани, параллельные $F$. И наоборот, каждая $(d-2)$-мерная грань
$Q$, параллельная $F$, содержится ровно в двух гипергранях пояска.
\end{rem}

\begin{sta}[Минковский \cite{Min1897}, Венков \cite{Ven1954}]
Выпуклый многогранник $P$ является параллелоэдром тогда и только тогда, когда он удовлетворяет условиям 1 -- 3.
\begin{enumerate}
	\item $P$ имеет центр симметрии.
	\item Любая гипергрань $P$ имеет центр симметрии.
	\item Любой поясок $P$ состоит из четырех или из шести гиперграней.
\end{enumerate}
\end{sta}

Основной результат данной статьи --- Теорема \ref{ordred}, предоставляющая комбинаторный критерий представимости параллелоэдра в виде прямого произведения
параллелоэдров меньшей размерности.

\section{Определения и обозначения}

\begin{defin}\label{antipod}
{\it Антиподальными} будем называть две гиперграни параллелоэдра $P$, если они переходят друг в друга при центральной симметрии $P$.
\end{defin} 

\begin{defin}\label{ven}
{\it Графом Венкова} параллелоэдра $P$ назовем граф $G_P$, если
\begin{enumerate}

\item Вершины $G_P$ находятся во взаимно-однозначном соответствии с (неупорядоченными) парами антиподальных гиперграней $P$.

\item Вершина $v_1$, соответствующая паре гиперграней $(F_1, F'_1)$, соединена красным ребром с вершиной $v_2$, соответствующей паре гиперграней
$(F_2, F'_2)$ тогда и только тогда, когда грани $F_1, F'_1, F_2, F'_2$ входят в один и тот же 6-поясок.

\item Вершина $v_1$, соответствующая паре гиперграней $(F_1, F'_1)$, соединена синим ребром с вершиной $v_2$, соответствующей паре гиперграней
$(F_2, F'_2)$ тогда и только тогда, когда грани $F_1, F'_1, F_2, F'_2$ образуют 4-поясок.
\end{enumerate}

\end{defin}

\begin{defin}\label{vencol}

Под {\it синим} или {\it красным} графом Венкова будем понимать подграф графа Венкова на том же множестве вершин, состоящий из всех синих или всех
красных ребер графа Венкова соответственно. Обозначим эти графы $G_P^b$ и $G_P^r$ соответственно.

\end{defin}

\begin{thm}[А.~Ордин \cite{Ord2005}]\label{ordred}

Параллелоэдр $P$ представляется в виде прямого произведения параллелоэдров меньшей размерности $P = P_1 \times P_2$ тогда и только тогда,
когда $G_P^r$ несвязен.

\end{thm}

\section{Граф Венкова произведения параллелоэдров.}

\begin{lem}\label{leftwards}

Пусть параллелоэдр $P$ представлен в виде прямого произведения двух параллелоэдров меньшей размерности: $P = P_1 \times P_2$. Тогда $G_P$
получается дизъюнктным объединением графов $G_{P_1}$ и $G_{P_2}$ и последующим проведением всех синих ребер $(v_1, v_2)$ таких, что
$v_1\in v(G_{P_1})$, $v_2\in v(G_{P_2})$.

\end{lem}

\begin{proof} 

Пусть $dim\, P = d$, $dim\, P_1 = d_1$, $dim\, P_2 = d_2$.

Любую гипергрань $F^{d-1}$ параллелоэдра $P$ можно представить либо как $F_1^{d_1-1} \times P_2$, либо как $P_1 \times F_2^{d_2-1}$, т.е.
как произведение гиперграни одного из параллелоэдров, входящих в разложение $P$, на другой параллелоэдр. Антиподальной $F^{d-1}$ гранью
будет, соответственно,  $(F')^{d-1} = (F'_1)^{d_1-1} \times P_2$ или $(F')^{d-1} = P_1 \times (F'_2)^{d_2-1}$, т.е. множитель, являющийся
гипергранью, заменяется на антиподальный. Таким образом, имеется естественное соответствие между $v(G_P)$ и 
$v(G_{P_1})\sqcup v(G_{P_2})$.

Опишем все пояски параллелоэдра $P$. Легко видеть, что любая его $(d-2)$-мерная грань $F^{d-2}$ однозначно задает поясок, в который входят обе гиперграни
$P$, содержащие $F^{d-2}$. С другой стороны, каждый поясок можно задать таким образом любой из $(d-2)$-мерных граней, по которой смежны соседние грани
этого пояска.

Для $(d-2)$-мерной грани $F^{d-2}$, задающей поясок, имеются 3 возможности. Разберем их отдельно.

\noindent {\bf Случай 1. } $F^{d-2} = F_1^{d_1-1} \times F_2^{d_2-1}$, т.е. $F^{d-2}$ является произведением гиперграней $P_1$ и $P_2$. Тогда соответствующий
поясок состоит из граней $F_1^{d_1-1} \times P_2$, $(F'_1)^{d_1-1} \times P_2$, $P_1 \times F_2^{d_2-1}$, $P_1 \times (F'_2)^{d_2-1}$,
где $(F'_1)^{d_1-1}$ и $(F'_2)^{d_2-1}$ --- гиперграни параллелоэдров $P_1$ и $P_2$ соответственно, антиподальные $F_1^{d_1-1}$ и $F_2^{d_2-1}$.
Такой поясок соответствует одному синему ребру между $G_{P_1}$ и $G_{P_2}$, причем, как легко видеть, проведены все возможные ребра.

\noindent {\bf Случай 2. } $F^{d-2} = F_1^{d_1-2} \times P_2$. Рассмотрим поясок $P_1$, заданный гранью $F_1^{d_1-2}$ коразмерности 2. Легко видеть,
что, умножив все гиперграни этого пояска на $P_2$, получим поясок $P$, соответствующий $F^{d-2}$. Этим пояскам соответствуют в точности
ребра $G_{P_1}$.

\noindent {\bf Случай 3. } $F^{d-2} = P_1 \times F_2^{d_2-2}$. Аналогично случаю 2, пояскам, заданным такими гранями, соответствуют в точности
ребра $G_{P_2}$.

Таким образом показано, что $G_P$ имеет структуру, описанную в условии Леммы \ref{leftwards}.

\end{proof}
 
\section{Фацетные векторы параллелоэдра, имеющего несвязный $G_P^r$.}

\begin{defin}\label{facvec}

Пусть $F$ --- гипергрань параллелоэдра $P$, $P'$ --- параллелоэдр, смежный с $P$ по $F$. Вектор $\mathbf t(F)$, удовлетворяющий
соотношению
$$P' = P + \mathbf t(F),$$
называется {\it фацетным вектором} грани $F$.

\end{defin}

Таким образом, каждой вершине $G_P$ можно сопоставить пару противоположных векторов --- фацетные векторы соответствующей пары гиперграней.

Обозначим через $\Lambda$ решетку, порожденную всеми фацетными векторами $P$, т.е. решетку центров всех параллелоэдров разбиения. 
Пусть $A \subset v(G_P)$. Символом $\Lambda(A)$ будем обозначать подрешетку решетки $\Lambda$, порожденную всеми фацетными векторами,
соответствующими вершинам множества $A$.

\begin{lem}\label{dirsum}

Пусть $v(G_P) = A_1\sqcup A_2$ и между $A_1$ и $A_2$ нет ни одного красного ребра. Тогда
$$\Lambda = \Lambda(A_1) \oplus \Lambda(A_2).$$

\end{lem}

\begin{proof}

Поскольку $\Lambda(A_1) + \Lambda(A_2)$ --- это решетка, порожденная всеми фацетными векторами, очевидно,
$$\Lambda = \Lambda(A_1) + \Lambda(A_2).$$

Таким образом, остается доказать, что $\Lambda(A_1) \cap \Lambda(A_2) = \{\mathbf 0 \}$.

Для доказательства потребуется ввести новое понятие.

\begin{defin}

Пусть имеется локально конечное нормальное полиэдральное разбиение $T$ пространства $\mathbb R^d$, и на его ячейках задана функция 
$f: \mathcal F^d (T) \to \mathbb R$. {\it Функцией приращения} функции $f$ называется такая функция ориентированных гиперграней ячеек
$g_f: \mathcal F^{d-1}_+ (T) \to \mathbb R$, что для любой ориентированной гиперграни $F^{d-1}_+$ выполняется
$$g_f(F^{d-1}_+) = f(P_2) - f(P_1),$$
где $P_1$ и $P_2$ --- ячейки $T$, смежные по этой гиперграни, и переход от $P_1$ к $P_2$ через $F^{d-1}_+$ происходит в положительном направлении.
В этом случае вместо записи $g_f(F^{d-1}_+)$ будем также иногда использовать запись $g_f(P_1, P_2)$.

\end{defin}

Легко видеть, что две функции на ячейках разбиения имеют одинаковые функции приращения тогда и только тогда, когда они различаются на константу.

Известно следующее утверждение.

\begin{sta}[Рышков и Рыбников \cite{RRy1997}]\label{rr}

Функция $g: \mathcal F^{d-1}_+ (T) \to \mathbb R$ является функцией приращения для некоторой функции $f: \mathcal F^d (T) \to \mathbb R$ тогда и только
тогда, когда 
\begin{enumerate}
\item[\rm 1.] $g(P_1, P_2) = -g(P_2, P_1)$ для любой пары смежных по гиперграни ячеек $P_1, P_2$.
\item[\rm 2.] Для любой $(d-2)$-мерной грани разбиения $F^{d-2}$ выполнено соотношение
$$g(P_1, P_2) + g(P_2, P_3) + \ldots + g(P_m, P_1) = 0,$$
где $P_1, P_2, \ldots, P_m$ --- все ячейки разбиения, содержащие $F^{d-2}$, и каждая из этих ячеек смежна по гиперграни с двумя соседними по циклу.
\end{enumerate}
\end{sta} 

Теперь приступим к доказательству леммы \ref{dirsum}. Предположим, что найдется такая точка $\lambda \in \Lambda \setminus \{\mathbf 0 \}$, что
$$\lambda \in \Lambda(A_1) \cap \Lambda(A_2).$$

Рассмотрим параллелоэдры разбиения $P_1$ и $P_2 = P_1 + \lambda$. Принадлежность $\lambda \in \Lambda(A_1)$ означает, что от $P_1$ к $P_2$
можно добраться, переходя только через гиперграни, соответствующие компоненте $A_1$. Аналогично, из $\lambda \in \Lambda(A_2)$ следует, что
между $P_1$ и $P_2$ существует путь, переходящий только через гиперграни, соответствующие компоненте $A_2$.

Рассмотрим какую-либо функцию $f$ на параллелоэдрах разбиения, сопоставляющую каждому параллелоэдру линейную функцию его центра (одну и ту же
для всех параллелоэдров). При этом линейную функцию выберем так, чтобы выполнялось $f(P_1) \neq f(P_2)$.

Выберем $\alpha\neq 1$ и построим функцию $g : \mathcal F^{d-1}_+ (T) \to \mathbb R$ следующим образом:
$$g(F^{d-1}_+) = \left\{ 
\begin{array}{ll}
g_f(F^{d-1}_+),\quad & \text{если $F^{d-1}$ соответствует компоненте $A_1$,}\\
\alpha g_f(F^{d-1}_+),\quad & \text{если $F^{d-1}$ соответствует компоненте $A_2$.}\\
\end{array}
\right.
$$

Легко проверить, что условия утверждения \ref{rr} для функции $g$ выполнены, следовательно, $g = g_{\phi}$ для некоторой функции 
$\phi : \mathcal F^d (T) \to \mathbb R$.

Рассмотрим путь от $P_1$ к $P_2$ через гиперграни, соответствующие $A_1$. Поскольку изменения функций $f$ и $\phi$ вдоль этого пути
одинаковы, имеем $\phi(P_2) - \phi(P_1) = f(P_2) - f(P_1)$. Аналогичным образом, рассматривая путь от $P_1$ к $P_2$ через гиперграни, соответствующие 
$A_2$, получим $\phi(P_2) - \phi(P_1) = \alpha(f(P_2) - f(P_1))$. Это возможно только при $f(P_2) - f(P_1) = 0$, что противоречит выбору $f$.

Таким образом, предположение о пересечении $\Lambda(A_1)$ и $\Lambda(A_2)$, отличном от нулевого, несостоятельно, и Лемма \ref{dirsum} доказана.

\end{proof}

\section{Невыпуклые множества с полиэдральной границей.}

Пусть $K\subset \mathbb R^d$. Через $cl(K)$ будем обозначать замыкание множества $K$, а через $int\, K$ --- множество его внутренних точек.

\begin{defin}\label{reg}

{\it Регулярным множеством} будем называть непустое замкнутое множество $K\subset \mathbb R^d$, если выполняется условие $cl(int\, K) = K$.

\end{defin}

\begin{lem}\label{triang}
Пусть дано невыпуклое регулярное множество $K$ с линейно-связной внутренностью. Тогда найдутся такие 3 точки $x,y,z \in int\, K$, 
что 
\begin{equation}\label{xyz}
[x,y]\subset int\, K, \quad [y,z] \subset int\, K, \quad [x,z] \not \subset K.
\end{equation}
\end{lem}

\begin{proof}

Покажем сначала, что найдется такая последовательность точек $x_0, x_1, \ldots, x_m$, что
\begin{equation}\label{broken}
[x_0, x_m] \notin K, \quad \text {и} \quad [x_{i-1}, x_i] \in int\, K \quad \text{при всех} \quad i=1,2\ldots m.
\end{equation}

Поскольку $K$ невыпукло, найдутся такие точки $x,z \in K$, что 
$$w = \alpha x + (1-\alpha) z \notin K \quad \text{для некоторого} \quad \alpha \in (0,1).$$

Из замкнутости $K$ следует, что $w$ не содержится в $K$ вместе с некоторой окрестностью. Кроме того, из условия $cl(int\, K) = K$ следует,
что существуют точки $x', y' \in int\, K$, сколь угодно близкие к $x$ и $y$ соответственно. Выберем $x'$ и $y'$ так, чтобы
$$\alpha x + (1-\alpha) z \notin K.$$

Из линейной связности $int\, K$ следует, что $x'$ и $y'$ можно соединить кривой, целиком лежащей в $int\, K$. Поскольку $int\, K$ открыто,
данную кривую можно заменить на некоторую вписанную в нее ломаную, также целиком лежащую в $int\, K$. Обозначив вершины этой ломаной в
порядке следования через $x_0 = x', x_1, \ldots, x_m = y'$, получим искомую последовательность.

Рассмотрим последовательность наименьшей длины, удовлетворяющую условиям (\ref{broken}). Покажем, что $m=2$.

Пусть $m>2$. Тогда имеется 2 возможных случая.

\noindent {\bf Случай 1. } $[x_0, x_{m-1}]\subset int\, K$. Тогда последовательность $x_0, x_{m-1}, x_m$ удовлетворяет условиям (\ref{broken}),
и имеет меньшую длину. Противоречие с выбором последовательности.

\noindent {\bf Случай 2. } $[x_0, x_{m-1}]\not\subset int\, K$. Конкретнее, существует такое $\alpha \in (0,1)$, что
$$w = \alpha x_0 + (1-\alpha)x_{m-1} \notin int\, K.$$

У точки $w$ существует сколь угодно близкая точка $w'\notin K$. Возьмем $w'$ так, чтобы для точки 
$$x'_{m-1} = x_{m-1}+\frac{1}{1-\alpha}(w' - w)$$
выполнялось $[x_{m-2}, x'_{m-1}]\subset int\, K$. Тогда нашлась последовательность $x_0$, $x_1,~\ldots$, $x_{m-2}$, $x'_{m-1}$, меньшая
исходной по длине и удовлетворяющая условиям (\ref{broken}). Снова получено противоречие с выбором исходной последовательности.

Таким образом, $m=2$, и можно положить $x=x_0, y=x_1, z=x_2$. 

\end{proof}

\begin{defin}\label{locconv}
Множество $K$ называется {\it локально выпуклым} в точке $x$, если существует такое открытое множество $U \supset x$, что $K \cap U$ выпукло. 
Если такого открытого множества не существует, $K$ невыпукло в точке $x$.
\end{defin}

\begin{defin}
Пусть $\mathcal K$ --- полиэдр. Объединение всех замкнутых $k$-мерных граней $\mathcal K$ будем называть $k$-ым {\it скелетом} $\mathcal K$
и обозначать $sk_k(\mathcal K)$.
\end{defin}

\begin{lem}\label{ridge}
Пусть невыпуклое регулярное множество $K \subset \mathbb R^d$ имеет линейно-связную внутренность, и $\partial K$ --- локально конечный 
полиэдр. Тогда существует точка 
$$x \in (sk_{d-2}(\partial K) \setminus sk_{d-3} (\partial K)),$$ 
в которой $K$ невыпукло.
\end{lem}

\begin{proof}

По лемме \ref{triang}, выберем точки $x,y,z \in int\, K$, удовлетворяющие условиям (\ref{xyz}). Легко видеть, что тем же самым условиям
удовлетворяет и любая достаточно близкая тройка точек $(x', y', z')$.

\begin{defin}\label{transv}

Аффинные пространства $L_1, L_2 \subset \mathbb R^d$ назовем {\it трансверсальными}, если:
\begin{enumerate}

\item $dim\, \aff (L_1\cup L_2) = \min ( dim\, L_1 + dim\, L_2 + 1, d).$

\item $dim\, (lin\, L_1 \cap lin\, L_2) = \max (d - dim\, L_1 - dim\, L_2, -1)$ в предположении, что $dim\, \varnothing = -1$.

\end{enumerate}

\end{defin}

Легко видеть, что почти любая (например, в смысле меры Лебега в $\mathbb R^{3d}$)
тройка точек в достаточно малой окрестности тройки $(x, y, z)$ задает двумерную плоскость, трансверсальную всем граням $\partial K$. Следовательно,
не умаляя общности, можно считать, что и $\aff \{x,y,z\}$ трансверсальна всем граням $\partial K$.

Легко видеть, что $conv\, \{x,y,z\} \cap \partial K \neq \varnothing$.

Среди всех точек множества $conv\, \{x,y,z\} \cap \partial K$ выберем точку $u$, наиболее удаленную от прямой $\aff \{x,z\}$, а если такая точка не
единственна, то из них возьмем наиболее близкую к прямой $\aff \{x,y\}$. Покажем, что $K$ локально невыпукло в точке $u$.

В самом деле, $u \notin [x,y] \cup [y,z]$, поскольку $[x,y] \cup [y,z] \subset int\, K$. Кроме того, $[x,z] \setminus K \neq \varnothing$.
Значит, в силу выбора $u$ имеем $d(u, [x,z]) \neq 0$. Следовательно, 
$$u \in int\, conv\, \{x,y,z\}.$$

Предположим теперь, что множество $K \cap conv\, \{x,y,z\}$ локально выпукло в точке $u$. Это означает, что существует прямая $\ell \in \aff (x,y,z)$,
являющаяся опорной в точке $u$ к множеству 
$$K \cap conv\, \{x,y,z\} \cap U,$$
где $U$ --- некоторая открытая окрестность точки $u$. 

С другой стороны, из выбора точки $u$ следует, что
$$\bigl\{u'\in conv\, \{x,y,z\}: d(u', [x,z]) > d(u, [x,z]) \bigr \} \subset K.$$
Поэтому $\ell$ может быть опорной прямой только в случае $\ell \parallel [x,z]$. Но тогда, еще раз пользуясь выбором точки $u$, получаем
$$\ell \cap U \cap conv\, \{x,y,z\} \subset \partial K.$$

Последнее включение противоречит выбору точки $u$, поскольку на прямой $\ell$ есть точки $\partial K$, более близкие к $[x,y]$, чем $u$. Полученное противоречие
показывает, что $K \cap conv\, \{x,y,z\}$ локально выпукло в точке $u$. Поэтому $K$ также невыпукло в точке $u$.

Далее заметим, что $u\notin (sk_{d-1}(\partial K) \setminus sk_{d-2}(\partial K))$, поскольку $K$ является локально выпуклым во всех внутренних
точках своих гиперграней. Кроме того, $u\notin sk_{d-3}(\partial K)$ в силу трансверсальности $\aff \{x,y,z\}$ всем $(d-3)$-мерным граням
$\partial K$, что как следствие влечет $sk_{d-3}(\partial K) \cap \aff \{x,y,z\} = \varnothing$.

Таким образом, $u \in (sk_{d-2}(\partial K) \setminus sk_{d-3} (\partial K))$, и лемма \ref{ridge} доказана.

\end{proof}

\begin{rem}
Теорема Бердона \cite[Theorem 7.5.1]{Bea1995} утверждает, что множество $K\subset \mathbb R^d$ выпукло тогда и только тогда, когда $K$ локально 
выпукло и связно. Лемма \ref{ridge} является уточнением этой теоремы в случае полиэдрального регулярного множества $K$.
\end{rem}

\section{Выделение прямых слагаемых для параллелоэдра с несвязным $G_P^r$.}

\begin{lem}\label{product}

Пусть $v(G_P) = A_1\sqcup A_2$ и между $A_1$ и $A_2$ нет ни одного красного ребра. Тогда

\begin{enumerate}

\item [\rm 1.] Множество $P+\Lambda(A_1)$ выпукло.

\item [\rm 2.] Имеет место представление в виде прямой суммы Минковского
$$P+\Lambda(A_1) = \aff \Lambda(A_1) \oplus P_2,$$
где $P_2 \subset \aff \Lambda(A_2)$ --- некоторый параллелоэдр.
\end{enumerate}

\end{lem}

\begin{proof}

В справедливости пункта 1 легко убедиться прямой проверкой условий леммы \ref{ridge} для множества $K = P+\Lambda(A_1)$. Действительно, возможны 2
случая невыпуклости $P+\Lambda(A_1)$ в $(d-2)$-мерной грани $F^{d-2}$.

\noindent {\bf Случай 1.} $F^{d-2}$ соответствует 6-пояску $P$, причем
$$\{\lambda : \lambda \in \Lambda(A_1),\; F^{d-2} \subset P+\lambda \} = \{\lambda_1, \lambda_2\},$$
$$\{\lambda : \lambda \in \Lambda,\; F^{d-2} \subset P+\lambda \} = \{\lambda_1, \lambda_2, \lambda_3\}.$$

Тогда для любых $i,j = 1,2,3$, $i\neq j$ вектор $\lambda_i - \lambda_j$ является фасетным вектором $P$. Имеем
$\lambda_3 \notin \Lambda(A_1)$, а $\lambda_1, \lambda_2 \in \Lambda(A_1)$, следовательно
$$\lambda_3 - \lambda_1 \in \Lambda(A_2) \qquad \text{и} \qquad \lambda_3 - \lambda_2 \in \Lambda(A_2).$$
Таким образом, $\lambda_2 - \lambda_1 \in \Lambda(A_2)$. С другой стороны, $\lambda_2 - \lambda_1 \in \Lambda(A_1)$, что противоречит Лемме \ref{dirsum}.
Случай 1 невозможен. 

\noindent {\bf Случай 2.} $F^{d-2}$ соответствует 4-пояску $P$, причем
$$\{\lambda : \lambda \in \Lambda(A_1),\; F^{d-2} \subset P+\lambda \} = \{\lambda_1, \lambda_2, \lambda_3\},$$
$$\{\lambda : \lambda \in \Lambda,\; F^{d-2} \subset P+\lambda \} = \{\lambda_1, \lambda_2, \lambda_3, \lambda_4\},$$
и $P+\lambda_1$ имеет общую гипергрань с $P+\lambda_2$ и $P+\lambda_3$.

Тогда в силу аргументов, примененных в случае 1, $\lambda_2 - \lambda_1 \in \Lambda(A_1)$ и $\lambda_4 - \lambda_3 \in \Lambda(A_2)$.
Но $\lambda_2 - \lambda_1 = \lambda_4 - \lambda_3$. Получено противоречие с Леммой~\ref{dirsum}, и случай 2 также невозможен.

Докажем пункт 2. Пусть $\pi : \mathbb R^d \to \aff \Lambda(A_2)$ --- параллельная проекция в направлении $\aff \Lambda(A_1)$. Покажем справедливость
равенства
\begin{equation}\label{ds1}
 P+\Lambda(A_1) = \aff \Lambda(A_1) \oplus \pi (P). 
\end{equation}

Включение 
$$ P+\Lambda(A_1) \subset \aff \Lambda(A_1) \oplus \pi (P) $$
очевидно. Убедимся в обратном включении. Предположим, что нашлась такая точка $x+\mu \notin P + \Lambda(A_1)$, что $x\in \pi(P)$, 
$\mu\in \aff \Lambda(A_1)$.

Поскольку $x\in \pi(P)$, существует такое $\nu \in \aff \Lambda(A_1)$, что $x+\nu \in P$. Следовательно, для всех $\lambda \in \Lambda(A_1)$
выполняется
$$x+\nu+\lambda \in P+\Lambda(A_1).$$

В силу выпуклости множества $P+\Lambda(A_1)$ имеем
$$conv\, \{x+\nu+\lambda \mid \lambda\in \Lambda(A_1)\} \subset P+\Lambda(A_1).$$

Очевидно, что
$$conv\, \{x+\nu+\lambda \mid \lambda\in \Lambda(A_1)\} = x + \aff \Lambda(A_1) \ni x+\mu.$$

Таким образом, получено противоречие, и соотношение (\ref{ds1}) доказано.

Для завершения доказательства леммы \ref{product} осталось показать, что проекция $\pi(P)$~--- параллелоэдр. Действительно, семейство
$$\{P+\Lambda(A_1)+\lambda \mid \lambda\in \Lambda(A_2)\}$$
образует разбиение $\mathbb R^d$, причем допускающее факторизацию по $lin\, \Lambda(A_1)$. Факторразбиение дает разбиение
$\aff \Lambda(A_2)$ на параллельные копии $\pi(P)$. Таким образом, $P_2 = \pi(P)$ --- параллелоэдр.

\end{proof}

\begin{proof}[Доказательство Теоремы \ref{ordred}] Имеем:
\begin{multline*}
P = (P+\Lambda(A_1))\cap (P+\Lambda(A_2)) = \\(\aff \Lambda(A_1) \oplus P_2) \cap (P_1 \oplus \aff \Lambda(A_2)) = 
P_1 \oplus P_2.
\end{multline*}

В самом деле, первое равенство следует из леммы \ref{dirsum}, второе --- из леммы \ref{product}, наконец, третье равенство следует из определения прямой
суммы. 

Но прямая сумма Минковского и прямое произведение параллелоэдров суть одно и то же. Поэтому $P = P_1 \times P_2$.

\end{proof}

\end{document}